# A Graph-Theoretical Approach to Ring Analysis: An Exploration of Dominant Metric Dimension in Compressed Zero Divisor Graphs and Its Interplay with Ring Structures

Nasir Ali[a,*], Hafiz Muhammad Afzal Siddiqui[a], Muhammad Imran Qureshi[b]
[a]Department of Mathematics, COMSATS University Islamabad, Lahore Campus, Pakistan.
[b]Department of Mathematics, COMSATS University Islamabad, Vehari Campus, Pakistan.
**Email address:**
nasirzawar@gmail.com  (Nasir Ali),
hmasiddiqui@gmail.com (Hafiz Muhammad Afzal Siddiqui),
imranqureshi18@gmail.com  (Muhammad Imran Qureshi),
**Email address and ORCID ID's:**
https://orcid.org/0000-0003-4116-9673  (Nasir Ali)
https://orcid.org/0000-0003-1794-6460  (Hafiz Muhammad Afzal Siddiqui)
https://orcid.org/0000-0002-0681-6313  (Muhammad Imran Qureshi)

**Abstract:**

The paper systematically classifies rings based on the dominant metric dimensions (Ddim) of their associated CZDG, establishing consequential bounds for the Ddim of these compressed zero-divisor graphs. The authors investigate the interplay between the ring-theoretic properties of a ring $R$ and associated CZDG. An undirected graph consisting of vertex set $Z(R_E)\setminus\{[0]\} = R_E\setminus\{[0],[1]\}$, where $R_E = \{[x] : x \in R\}$ and $[x] = \{y \in R : ann(x) = ann(y)\}$ is called a compressed zero-divisor graph, denoted by $\Gamma_E(R)$. An edge is formed between two vertices $[x]$ and $[y]$ of $Z(R_E)$ if and only if $[x][y] = [xy] = [0]$, that is, iff $xy = 0$. For a ring $R$, graph $G$ is said to be realizable as $\Gamma_E(R)$ if $G$ is isomorphic to $\Gamma_E(R)$. Moreover, an exploration into the Ddim of realizable graphs for rings is conducted, complemented by illustrative examples reinforcing the presented results. A recent discussion within the paper elucidates the nuanced relationship between Ddim, diameter, and girth within the domain of compressed zero-divisor graphs. This research offers a comprehensive and insightful analysis at the intersection of algebraic structures and graph theory, providing valuable contributions to the current mathematical discourse.

**Keywords:** dominant metric dimensions, domination, compressed zero-divisor graph, resolvability, Algebraic structures

## 1. Introduction

Graph Exploring the intricate interplay between graph theory and algebraic structures has been a cornerstone of contemporary mathematics. The amalgamation of these two domains has not only enriched theoretical frameworks but also found practical applications across various disciplines. From modeling complex interactions in social networks to optimizing transportation routes and understanding molecular structures in chemistry, the versatility of graph theory knows no bounds [13]. Moreover, its intersection with algebra has led to the emergence of algebraic graph theory, a field that delves into the profound relationships between graphs and algebraic structures such as groups, matrices, and rings [3]. Among the myriad of algebraic structures, rings stand out as fundamental objects of study. Rings encapsulate essential algebraic properties, offering a rich tapestry for exploration. In recent years, there has been a growing interest in characterizing rings based on certain graph-theoretic properties. Specifically, the focus has been on elucidating the connection between the algebraic properties of a ring and the structural attributes of its associated graph. One such graph of paramount importance is the compressed zero-divisor graph (CZDG), denoted as $\Gamma_E(R)$. The CZDG of a ring R is an undirected graph that captures significant information about the zero-divisors of R in a compressed manner. Each vertex of the CZDG represents an equivalence class of zero-divisors, offering a compact representation of the zero-divisor structure of the ring. The edges of the CZDG encode the relationships between these equivalence classes, shedding light on the underlying algebraic properties of the ring.

This paper aims to address the following fundamental problem: Given a ring $R$, how can one characterize rings based on dominant metric dimensions (Ddim) of its associated compressed zero-divisor graph, denoted as $\Gamma_E(R)$? The investigation extends beyond the mere computation of Ddim, delving into the intricate interplay between the ring-theoretic properties of $R$ and the structural attributes of its CZDG. Furthermore, the paper seeks to classify rings based on the Ddim of their CZDGs, establishing bounds for these dimensions and providing insights into the realizable graphs of rings in the context of $\Gamma_E(R)$. The exploration of these aspects contributes to a deeper comprehension of the underlying relationships between algebraic structures and graph-theoretic representations.

Building upon seminal works in the field, such as those by Beck [7], Anderson, Livingston [3], Redmond, and Mulay [22], we embark on a journey to unravel the intricate connections between rings and their associated CZDGs. Through a careful analysis of these structures, we aim to uncover new insights into the underlying relationships between algebra and graph theory. The interplay between graph theory and algebra has been a subject of significant interest, with various scholars contributing to its exploration. One notable contribution is attributed to Beck, who introduced the concept of a zero-divisor graph (ZD-graph) for commutative rings. Beck's work, highlighted in [7], delved into the relationship between graph theory and algebra, particularly focusing on the correspondence between ring elements and graph nodes. In this framework, a zero vertex is intricately linked to all other vertices, presenting a unique perspective on graph coloring. Building upon Beck's foundation, Anderson and Livingston, as discussed in [3], extended the study of ZD-graphs to encompass nonzero zero divisors within commutative rings. Their investigation led to the development of an undirected graph, denoted Γ(R), where vertices represent nonzero zero divisors, and edges connect vertices x and y if xy equals zero. However, the exploration of ZD-graphs has evolved beyond its initial formulations. Recent advancements, as evidenced in [15] and [16], have broadened the scope of ZD-graphs to include novel variations such as ideal-based and module-based ZD-graphs. These innovations have revitalized the field, offering fresh perspectives and avenues for exploration. Moreover, Redmond's pioneering work, detailed in [22], expanded the application of ZD-graphs from unital commutative rings to noncommutative rings. Redmond introduced diverse methodologies to characterize ZD-graphs associated with noncommutative rings, encompassing both undirected and directed graph structures. Notably, Redmond extended this paradigm by introducing ideal-based ZD graphs, aiming to generalize the methodology by incorporating elements with zero products belonging to specific ideals of the ring. Mulay [16] work inspired us to study the ZD-graph obtained by considering equivalence classes of ZD of a ring $R$. This typed of ZD-graph is called CZDG, denoted by $\Gamma_E(R)$ [4]. A CZDG is an undirected graph obtained by considering $Z(R_E)\setminus\{[0]\} = R_E\setminus\{[0],[1]\}$ as vertex set, and can be constructed by taking the equivalence classes $[x] = \{y \in R : ann(x) = ann(y)\}$, for every $x \in R\setminus([0]\cup[1])$ as vertices and and edge is formed between two distinct classes $[x]$ and $[y]$ iff $[x][y] = 0$, i.e., iff $xy = 0$. It is important to note that if two vertices say $x$ and $y$ are adjacent in $\Gamma(R)$, then in CZDG, $[x]$ and $[y]$ are adjacent iff $[x] \neq [y]$. Clearly, $[1] = R\setminus Z(R)$ and $[0] = \{0\}$, also for each $x \in R \setminus ([0] \cup [1])$, $[x] \subseteq Z(R)\setminus\{0\}$. Readers may study [5] for some interesting results on CZDG. We consider an example to understand concept of ZD-graph and CZDG. Let $R = \mathbb{Z}_{16}$, then the vertex set of $\Gamma(R) =$

$\{2, 4, 6, 8, 10, 12, 14\}$, see Figure 1.1 (i) shows its ZD-graph. Now, we see $ann(2) = \{8\}$, $ann(4) = \{4, 8, 12\}$, $ann(6) = \{8\}$, $ann(8) = \{2, 4, 6, 8, 10, 12, 14\}$, $ann(10) = \{8\}$, $ann(12) = \{4,8,12\}$ and $ann(14) = \{6,8\}$. Hence, vertex set for $\Gamma_E(R) = \{[2], [4], [8], [14]\}$. see Figure 1.1 (ii) for its CZDG.

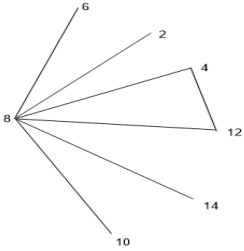
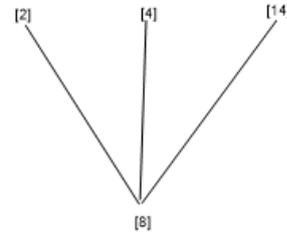

Figure 1.1 (i)     Figure 1.1 (ii)

The Annihilator ideals in the ring corresponds to vertices of the $\Gamma_E(R)$. Also remember that $diam(\Gamma_E(R)) \leq 2$ and CZDG is always connected. Also, $diam(\Gamma_E(R)) \leq diam(\Gamma(R))$. For the CZDG, $gr(\Gamma_E(R)) \leq 3$ whenever, CZDG of $R$ contains a cycle [5]. It can be seen in [4] that $\Gamma_E(R)=\Gamma_E(S)$ if $S$ is a Noetherian or $S$ is finite commutative ring. Readers may see [1, 2, 3] to read many advantages of studying CZDG over the earlier studied ZD-graph. For example, any ring $R$ having at least 2 vertices, there exists no finite regular CZDG [[27], Proposition 1.10]. Further, Spiroff et.al, [27] showed that the CZDG of local ring $R$ is isomorphic to a star graph with minimum 4 vertices. (If a ring $R$ has a unique maximal ideal then it is called local ring). Another crucial aspect to consider is the connection between studying equivalence classes's graphs. The associated primes of $R$ usually considered as distinct vertices in CZDG. In this paper, all graphs are simple graphs, a commutative ring with unity is denoted by $R$ and units set is considered as $U(R)$. $\mathbb{Z}_n$ denotes ring of integers modulo $n$ and $\mathbb{F}_q$ denotes finite field on $q$ elements. Readers are encouraged to study [11, 17] for basic definitions of graph theory and [6, 13] to study basic definitions of ring theory.

The graph associated with $R$ provides a remarkable demonstration of the properties of ZD the of lattice of ideals of $R$. This graph enables us to visualize and analyze the algebraic properties of rings through graph theory techniques. In [3], the properties of $\Gamma(R)$ were analyzed by the authors and found to be interesting. To study $\Gamma(R)$, we stick to the method presented by Anderson and Livingston in [3], where non-zero ZD are vertex set for $\Gamma(R)$. Unless otherwise stated, we consider $R$ is a finite unital commutative ring in this paper. $L(R)$ denotes set of non-zero ZD. Consider the ring $R$ contains only single maximal ideal then it is called local. In the context of a graph $G$, a subset $S \subseteq V(G)$ is defined as a dominating set (DS) if, for every vertex $x \in V(G)$ not in $S$, there exists at least one vertex $u \in S$ such that $x$ is adjacent to $u$. The smallest cardinality among all dominating sets of $G$ is known as the dominating number (DN) of $G$, denoted by $\gamma(G)$ [10]. The notion of metric dimension for connected graphs and its associated properties were initially introduced by Slater in 1975, followed independently by Harary and Melter in 1976. In the realm of graph theory, a subset of vertices denoted as S is deemed to resolve a graph $G$ if each vertex within $G$ can be uniquely determined by its vector of distances to the vertices in S. The smallest cardinality of such a resolving set for a given graph $G$ is termed a minimum resolving set, often referred to as a basis for $G$. The cardinality of this basis is formally recognized as the metric dimension of $G$, denoted by $dim(G)$. Mathematically, an ordered set $W = \{w_1, w_2, \ldots, w_k\} \subseteq V(G)$. If, for every pair of vertices $u$ and $v$ in $V(G)$, there is a unique representation with respect to $W$, i.e., $r(u|W) \neq r(v|W)$, where $r(u|W) = \big(d(u, w_1), d(u, w_2), \ldots, d(u, w_k)\big)$, then $W$ is referred to as a resolving set of $G$. The resolving set with the minimum cardinality is termed the metric dimension (MD) of G, denoted as $\dim(G)$. Brigham, et al. collectively examined MD and DS, coining the term "resolving domination number," denoted by $\gamma_r(G)$. They established that $\max\{\dim(G), \gamma(G)\} \leq \gamma_r(G) \leq \dim(G) + \gamma(G)$. In the realm of dominant sets, a dominant resolving set of G is an ordered set $W \subseteq V(G)$ such that $W$ serves as both a resolving set and a dominating set of $G$. The minimum cardinality of a dominant resolving set is termed a dominant basis of $G$, while the cardinality of the dominant basis is recognized as a dominant MD of G, denoted by $Ddim(G)$.

This research presents a pioneering exploration into the realm of dominant metric dimensions within the context of compressed zero-divisor graphs associated with rings. The novelty lies in the comprehensive analysis of the interplay between the algebraic properties of a ring and the corresponding compressed zero-divisor graph. The paper not only classifies rings based on the Ddim of their CZDGs but also provides bounds for these dimensions, shedding light on the inherent constraints governing the metric structure. Additionally, the investigation extends to the realm of realizable graphs, elucidating the conditions under which a graph is isomorphic to the compressed zero-divisor graph of a given ring. By doing so, the paper contributes valuable insights into the broader understanding of algebraic structures, offering a nuanced perspective on the relationship between dominant metric dimensions, diameter, and girth within the domain of compressed zero-divisor graphs. The results presented in this work pave the way for further exploration and applications in ring theory and graph theory, showcasing the significance of this study in advancing the current state of knowledge in mathematical research. Susilowati et al. [15] determined the Ddim of a specific class of graphs, characterized graphs with particular Ddim, and computed the Ddim of joint and comb products of graphs. In this context, we examine selected findings from references [11] as follows:

**Remark 1** [11]. In the context of path graphs, represented as $P_n$, and cyclic graphs, denoted as $C_n$, the dominating number is characterized by $\gamma(P_n) = \gamma(C_n) = \left\lceil \frac{n}{3} \right\rceil$. Additionally, the metric dimension for path graphs is given by $dim(P_n) = 1$, while for cyclic graphs, it is expressed as $dim(C_n) = 2$.

**Remark 2** [11]. In the context of complete graphs, represented as $K_n$, the dominating number is characterized by $\gamma(K_n) = 1$. Additionally, the metric dimension for complete graphs is given by $dim(K_n) = n - 1$.

**Remark 3** [11]. In the context of start graphs, represented as $S_n$, the dominating number is characterized by $\gamma(S_n) = 1$. Additionally, the metric dimension for star graphs is given by $dim(S_n) = n - 2, \forall\, n \geq 2$.

**Remark 4** [11]. In the context of complete bipartite graph, represented as $K_{m,n}$, the dominating number is characterized by $\gamma(K_{m,n}) = 2$. Additionally, the metric dimension for star graphs is given by $dim(K_{m,n}) = m + n - 2, \forall\, m, n \geq 2$.

Furthermore, we consider some previous results on dominant metric dimension of $G$ as follows:

**Theorem 1** [15]: For a cyclic graph $C_n$ with an order of $n \geq 7$, it holds that the dominating metric dimension of $C_n$ is equal to the dominating number, denoted as $\gamma(C_n)$.

**Theorem 2** [15]: For a star graph $S_n$ with an order of $n \geq 2$, the dominating metric dimension of the graph is determined to be $n - 1$.

**Theorem 3** [15]: Consider $K_{m,n}$ as a complete bipartite graph with $m$ and $n$ both greater than or equal to 2. In this case, the dominating metric dimension of $K_{m,n}$ is equivalent to its metric dimension, expressed as $Dim_d(K_{m,n}) = dim(K_{m,n})$.

**Theorem 4** [15]: For a path graph $P_n$, where n is greater than or equal to 4, the dominating metric dimension of the graph is equal to its dominating number, expressed as $Dim_d(P_n) = \gamma(P_n)$.

**Theorem 5** [15]: For a complete graph $K_n$, where $n$ is greater than or equal to 2, the dominating metric dimension of the graph is equal to its metric dimension, stated as $Dim_d(K_n) = dim(K_n)$.

**Theorem 6** [15]: The condition $Dim_d(P_n) = 1$ holds if and only if the graph $G$ is isomorphic to the path graph $P_n$, where n takes on values of 1 or 2.

Moreover, Ddim for a single vertex graph $G$ is supposed to be zero and for an empty graph it is undefined. Our discussion commences with the subsequent observation.

## 2. Dominant Metric dimension of some compressed zero-divisor graphs

**Proposition 2.1:** Let $R$ be a finite commutative ring. Then $Dim_d(\Gamma_E(R)) = 0$ iff $\Gamma(R) \cong K_n$.

**Proof.** Let $\Gamma(R) \cong K_n$, then one of the two cases holds, either $R$ is isomorphic to $\mathbb{Z}_2 \times \mathbb{Z}_2$ or for all $x, y \in Z^*(R)$ product of $x$ and $y$ is zero. i.e., $xy = 0$. Let $v_1, v_2, \ldots v_n$ be the ZD of the ring $R$ that correspond to the vertices of ZD-graph, then $[v_1] = [v_2], \cdots = [v_n]$ suggests that all vertices of ZD-graph would collapse to a graph with single vertex in CZDG, and as we know Ddim is 0 for a single vertex graph.

For the converse part, assume that $\Gamma(R) \not\cong K_n$. This suggests the presence of at least one vertex in $\Gamma(R)$ that is not connected to every other vertex. Consequently, the size of the edge set CZDG is at least 2, i.e., $|\Gamma_E(R)| \geq 2$, leading to $Dim_d(\Gamma_E(R)) \geq 1$.

The converse can also be shown by assuming that $Dim_d(\Gamma_E(R)) = 0$. This implies that $\Gamma_E(R) = \{[a]\}$ for some non-zero element say, $a \in Z^*(R)$, indicating that $\Gamma_E(R)$ is a single vertex graph. This implies that the ZD-graph $\Gamma(R)$ is either isomorphic to $K_n, \forall, n \geq 1$ or a graph with single vertex. □

**Proposition 2.2:** Let $R$ be a finite commutative ring. Then $Dim_d(\Gamma_E(R)) = 1$, iff $\Gamma(R) \cong K_{m,n}$, where $m$ or $n \geq 2$.

**Proof.** Suppose that the ZD-graph $\Gamma(R)$ is isomorphic to $K_{m,n}$ having two distance similar classes $V_1$ and $V_2$. Specifically, assume that $V_1 = \{u_1, u_2, \cdots, u_m\}$ and $V_2 = \{v_1, v_2, \cdots, v_m\}$ such that $u_i v_j = 0$, $\forall i \neq j$. Cleary, an independent set is formed by each of $V_1$ and $V_2$. Furthermore, observing that $[u_1] = [u_2] = \cdots = [u_m]$ and, $[v_1] = [v_2] = \cdots = [v_n]$, we deduce that both $V_1$ and $V_2$ each represent a single vertex CZDG. Given the connected nature of the graph, we can conclude that $\Gamma_E(R) \cong K_{1,1}$, which can be visualized as a path consisting of two vertices. Therefore, by Theorem 6, it follows that the Ddim of CZDG is 1 that is $Dim_d(\Gamma_E(R)) = 1$.□

**Remark 2.1:** It's important to note that the converse of Proposition 2.2 may not hold true; this is exemplified by the graph illustrated in Figure 2.1. However, in cases where $R \cong \mathbb{Z}_2 \times \mathbb{Z}_2$, we find that $\Gamma_E(R) \cong k_{1,1}$, with $Dim_d(\Gamma_E(R)) = 1$ such that $\Gamma(R) \cong \Gamma_E(R)$. As demonstrated in ([27], Proposition 1.5), a significant difference between ZD-graph and CZDG is that the latter cannot be isomorphic to a complete graph having minimum of 3 nodes. However, see ([27], Proposition 1.7), if $\Gamma_E(R)$ is isomorphic to a complete r-partite graph, then $r$ must equal to 2 resulting in $\Gamma_E(R) \cong k_{n,1}$, for some $n \geq 1$. A ring R is classified as a Boolean ring if $a^2 = a$ holds true for every element $a \in R$. Importantly, a Boolean ring $R$ is both commutative and has a characteristic of 2 (char$(R) = 2$). More comprehensively, a commutative ring qualifies as a von Neumann regular ring (VNR) when, for any element $a \in R$, there exists an element $b$ within $R$ such that $a = a^2 b$. This condition is equivalent to $R$ being a zero-dimensional reduced ring, as elucidated in ([13], Theorem 3.1). Clearly, we can classify Boolean ring as VNR; but, the converse may not always be true. For instance, consider a family $\{F_i\}_{i \in I}$ of fields, where the product $\prod_{i \in I} F_i$ is VNR. Nevertheless, it is Boolean if and only if $F_i \cong \mathbb{Z}_2$ holds for all $i \in I$. If $R$ is reduced ring, then for $r, s \in \Gamma(R)$, the conditions $N(r) = N(s)$ and $[r] = [s]$ are equivalent ([13], Lemma 3.1). Furthermore, if $R$ is a VNR, then these conditions are equivalent to $rR = sR$. Moreover, if $R$ is a VNR and $B(R)$ signifies the set of idempotent elements within $R$, the mapping defined by $e \mapsto [e]$ forms an isomorphism from the subgraph of $\Gamma(R)$ induced by $B(R)\setminus\{0,1\}$ onto $\Gamma_E(R)$ ([13], Proposition 4.5). Particularly, if $R$ is a Boolean ring (i.e., $R = B(R)$), then $\Gamma_E(R)$ is isomorphic to $\Gamma(R)$. This discourse leads us to the subsequent characterization.

**Corollary 2.1:** Consider $R$ & $S$ bring reduced commutative ring. If $\Gamma(R)$ is isomorphic to $\Gamma(S)$ then $Dim_d(\Gamma_E(R)) = Dim_d(\Gamma_E(S))$.

## 3. Analyzing bounds for the dominant metric dimension of $\Gamma_E(R)$

Here, we explore the importance of calculating Ddim in studying CZDG. Additionally, we determine the Ddim of certain specialized ring types that correspond to $\Gamma_E(R)$. Notably, a recent contribution by Pirzada et al. [18] presented a work on characterization of $\Gamma(R)$ in cases where the MD is finite and cases where it remains undefined ([18], Theorem 3.1).

**Proposition 3.1.** Let $R$ be a finite commutative ring with unity. Then
  i. $Dim_d(\Gamma_E(R))$ is finite iff $R$ is finite.
  ii. $Dim_d(\Gamma_E(R))$ is undefined iff $R$ is an integral domain.

**Proof.**
  i. Assume that $Dim_d(\Gamma_E(R))$ is finite, then $\exists$ a minimal dominant metric basis for $\Gamma_E(R)$, we denote this basis set as $\{v_1, v_2, \ldots, v_t\}$. By using ([2], Theorem 2.3), $diam(\Gamma_E(R)) \leq 3$. So, $d(r, e) = 0, 1, 2$ or $3$ for every $r \in V(\Gamma_E(R))$ and $e \in E(\Gamma_E(R))$. Hence, $|L(R)| \leq 4t$, which implies that $\Gamma_E(R)$ is finite, and hence $R$ is finite. Conversely, given that $R$ is finite, then $|L(R)|$ is finite, since $\Gamma_E(R)$ is contained in $R$. So, $Dim_d(\Gamma_E(R))$ is finite.
  ii. As we know if $R$ is an ID then $\Gamma_E(R)$ is not defined which follows that $Dim_d(\Gamma_E(R))$ is undefined and vice versa. ∎

Let us consider the following lemma,

**Lemma 3.1:** Consider a ring $R$ which is finite local ring, $|R| = p^n$, with $n$ being a positive integer and some prime $p$.

The preceding result will be applied to determine the Ddim of rings $R$ which are finite local rings. □

**Proposition 3.1:** Let $|R| = p^2$ for a local ring $R$ and $p$ being $2, 3$ or $5$, then $Dim_d(\Gamma_E(R))$ can be either undefined or $0$.

**Proof.** Let us assume all local rings having order $p^2$, where $p$ is a prime. The following rings $\mathbb{F}_{p^2}$, $\mathbb{Z}_{p^2}$ and $\frac{\mathbb{F}_p[x]}{(x^2)}$ are local rings of order $p^2$ ([10], p. 687).

Case I: If $R \cong \mathbb{F}_{p^2}$, i.e., $R$ is a field having order $p^2$. In such cases, the graph $\Gamma_E(R)$ becomes an empty graph, resulting in an undefined $Dim_d(\Gamma_E(R))$.

Case II: If $R \cong \mathbb{Z}_{p^2}$, or $R \cong \frac{\mathbb{F}_p[x]}{(x^2)}$, i.e., $|R| = p^2$ and $R$ is not a field. For $p = 2, 3$ or $5$ the graph $\Gamma_E(R)$ consists of a single vertex, leading to $Dim_d(\Gamma_E(R)) = 0$. This concludes our result. □

**Proposition 3.2:** Let $R$ be a local ring having order, (i.e., $R$ is not a field)

(i) $p^3$ with $p = 2$ or $3$, then $Dim_d(\Gamma_E(R)) = 0$, and $Dim_d(\Gamma_E(R)) = 1$ only if $R \cong \mathbb{Z}_8, \mathbb{Z}_{27}, \mathbb{Z}_2[x]/(x^3), \mathbb{Z}_3[x]/(x^3), \mathbb{Z}_4[x]/(2x, x^2 - 2), \mathbb{Z}_9[x]/(3x, x^2 - 3), \mathbb{Z}_9[x]/(3x, x^2 - 6)$.

(ii) $p^4$ with $p = 2$, then $Dim_d(\Gamma_E(R))$ is equal to $0, 1$ or finite.

**Proof.** (a) Following rings $\mathbb{F}_{p^3}$, $\frac{\mathbb{F}_p[x]}{(x^3)}$, $\frac{\mathbb{F}_p[x,y]}{(x,y)^2}$, $\frac{\mathbb{Z}_{p^2}[x]}{(px, x^2)}$, $\frac{\mathbb{Z}_{p^2}[x]}{(px, x^2-p)}$ are all local rings having order $p^3$.

Case (i). When $p = 2$, the local rings $\frac{\mathbb{Z}_2[x,y]}{(x,y)^2}$ and $\frac{\mathbb{Z}_4[x]}{(2x, x^2)}$ have same equivalence classes of zero divisors. Let us assume that $[m] = \{x, y, x + y\}$ is equivalence class for any ZD $m$ of ring $\frac{\mathbb{Z}_2[x,y]}{(x,y)^2}$ and $[n] = \{2, x, x + 2\}$ is equivalence class for any ZD $m$ ring $\frac{\mathbb{Z}_4[x]}{(2x, x^2)}$, i.e., $\Gamma_E(R)$ is a graph with single vertex. Hence, $Dim_d(\Gamma_E(R)) = 0$. But $\Gamma_E(R) \cong P_2$ for $\mathbb{Z}_8, \mathbb{Z}_2[x]/(x^3), \mathbb{Z}_4[x]/(2x, x^2 - 2)$. Hence by Theorem 6, $Dim_d(\Gamma_E(R)) = 1$.

Case (ii). If $p = 3$ for above given rings, we found that CZDG structures of the rings $\mathbb{Z}_{27}, \mathbb{Z}_3[x]/(x^3), \mathbb{Z}_9[x]/(3x, x^2 - 3), \mathbb{Z}_9[x]/(3x, x^2 - 6)$ is same and isomorphic to $P_2$. Then by Theorem 6, $Dim_d(\Gamma_E(R)) = 1$.

Also, the rings $\frac{\mathbb{Z}_{3^2}[x]}{(3x,x^2)}$ and $\frac{\mathbb{Z}_3[x,y]}{(x,y)^2}$ have same equivalence classes of ZD given by $[m] = \{3, 6, x, 2x, x + 3, x + 6, 2x + 3, 2x + 6\}$ for any ZD $m$ of ring $\frac{\mathbb{Z}_{3^2}[x]}{(3x,x^2)}$ and $[n] = \{x, 2x, y, 2y, x + y, 2x + y, x + 2y,$

$2x + 2y$} for any ZD $n$ of ring $\frac{\mathbb{Z}_3[x,y]}{(x,y)^2}$, i.e., $\Gamma_E(R)$ is a graph with single vertex. Hence, $Dim_d(\Gamma_E(R)) = 0$.

(ii) Now we focus our attention to local rings having order $p^4$, by taking $p = 2$. It is found in [10] that there are 21 non-isomorphic commutative local rings with identity of order 16. The following are rings with $Dim_d(\Gamma_E(R)) = 0$ are $\mathbb{F}_4[x]/(x^2)$, $\mathbb{Z}_2[x,y,z]/(x,y,z)^2$, and $\mathbb{Z}_4[x]/(x^2 + x + 1)$. The following are the rings having $Dim_d(\Gamma_E(R)) = 1$, $\mathbb{Z}_2[x]/(x^4)$, $\mathbb{Z}_2[x, y]/(x^3, xy, x^2)$, $\mathbb{Z}_4[x]/(2x, x^3 - 2)$, $\mathbb{Z}_4[x]/(x^2 - 2)$, $\mathbb{Z}_8[x]/(2x,\ x^2)$, $\mathbb{Z}_{16}$, $\mathbb{Z}_4[x]/(x^2 - 2x - 2)$, $\mathbb{Z}_8[x]/(2x,\ x^2 - 2)$ and $\mathbb{Z}_4[x]/(x^2 - 2x)$. Furthermore, the following rings $\mathbb{Z}_4[x]/(x^2)$, $\mathbb{Z}_2[x,y]/(x^2,y^2)$ and $\mathbb{Z}_2[x,y]/(x^2 - y^2,\ xy)$ have $Dim_d(\Gamma_E(R))$ is finite. □

Now we will find $Dim_d(\Gamma_E(\mathbb{Z}_n))$.

**Proposition 3.3:** Consider a prime number $p$.
(a) When $n = 2p$ & $p > 2$, then $Dim_d(\Gamma_E(\mathbb{Z}_n)) = 1$.
(b) When $n = p^2$, then $Dim_d(\Gamma_E(\mathbb{Z}_n)) = 0$.

**Proof.** (a) When $p = 2$, CZDG of $\mathbb{Z}_4$ is a single vertex graph, hence $Dim_d(\Gamma_E(\mathbb{Z}_n)) = 0$.

Consider $p > 2$, then $\{2, 2.2, 2.3, \ldots, 2.(p-1), p\}$ is ZD set of $\mathbb{Z}_n$. Given that the, $Char\ (\mathbb{Z}_n) = 2p$, we can deduce that $p$ is adjacent with every other vertex. Consequently, the equivalence classes of these ZD are as follows:

$$[2] = [2.2] = \ldots = [2.(p-1)] = \{p\}, \qquad [p] = \{2, 2.2, 2.3, \ldots, 2.(P-1)\}.$$

Thus, $\{[p], [2x]\}$ is vertex set of $\Gamma_E(\mathbb{Z}_n)$, for some positive integer $x = 1, 2, 3, \ldots, p - 1$. So $\Gamma_E(\mathbb{Z}_n) \cong P_2$ then, by Theorem 6, $Dim_d(\Gamma_E(\mathbb{Z}_n)) = 1$.

(b) The ZD set of $\mathbb{Z}_n$ is $\{p, p.2, p.3, \ldots, p.(p-1)\}$, when $n = p^2$ and $p > 2$. Since $Char\ (\mathbb{Z}_n) = p^2$, so equivalence classes of all are same, i.e., $\{p, p.2, p.3, \ldots, p.(p-1)\}$. Hence CZDG of $\mathbb{Z}_n$ is a single vertex graph, so $Dim_d(\Gamma_E(\mathbb{Z}_n)) = 0$.

Consider a ring $R$, a graph $G$ is said to be realizable as $\Gamma_E(R)$ if $G \cong \Gamma_E(R)$. However, numerous results suggest that the majority of graphs cannot be realized as $\Gamma_E(R)$, For instance, $\Gamma_E(R)$ cannot represent a complete graph with three or more vertices or cycle graph. □

**Proposition 3.4:** Consider a realizable graph (RG) $\Gamma_E(R)$ containing 3 vertices then $Dim_d(\Gamma_E(R)) = 1$.

**Proof.** In a study by Spiroff et al. [27], it was established that one RG, $\Gamma_E(R)$ containing precisely 3 vertices as a graph of equivalence classes of ZD for some ring $R$ is $P_3$, see Figure 3.1. understandably, $Dim_d(\Gamma_E(R)) = 1$. □

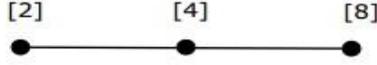

Figure 3.1: $\mathbb{Z}_{16}$

**Proposition 3.5:** Consider a RG $\Gamma_E(R)$ containing 4 vertices then $Dim_d(\Gamma_E(R))$ is either 1 or 2.

**Poof.** The possible RG $\Gamma_E(R)$ with four vertices are illustrated in Figure 3.2., it is readily apparent that their Ddim can be either 1 or 2. □

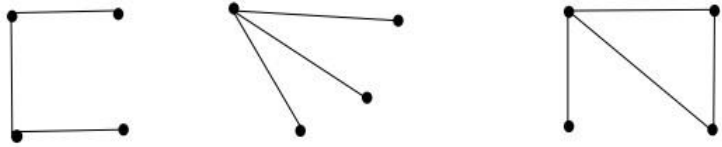

Figure 3.2: $(\mathbb{Z}_4 \times \mathbb{F}_4)$, $\mathbb{Z}_4[x]/(x^2)$, $\mathbb{Z}[x,y]/(x^3,xy)$

**Proposition 3.6:** Consider a RG $\Gamma_E(R)$ with 5 vertices then $Dim_d(\Gamma_E(R)) = 1$.

**Poof.** The possible RG $\Gamma_E(R)$ with 5 vertices are illustrated in Figure 3.3; it is readily apparent that their Ddim is 1. □

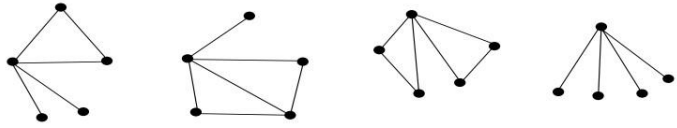

Figure 3.3: $\mathbb{Z}_9[x]/(x^2), \mathbb{Z}_{64}, \mathbb{Z}_3[x,y]/(xy,x^3,y^3,x^2-y^2), \mathbb{Z}_8[x,y]/(x^2,y^2,4x,\ 4y,\ 2xy)$

## 4. Relationship between dominant metric dimension, diameter and girth of $\Gamma_E(R)$

Within this section, we delve into the correlation among dominant metric dimension, diameter and girth of $\Gamma_E(R)$. Since $gr(\Gamma_E(R)) \in \{3, \infty\}$, but $gr(\Gamma_E(R)) = 3$ iff $gr(\Gamma(R)) = 3$, where $R$ is a reduced commutative ring. Moreover, $gr(\Gamma_E(R)) = \infty \Leftrightarrow gr(\Gamma(R)) \in \{4, \infty\}$. Nevertheless, it is possible that $gr(\Gamma(R)) = 3$ & we could have either $gr(\Gamma_E(R)) = 3$ or $\infty$. In the subsequent result, we establish the $Dim_d(\Gamma_E(R))$ of a ring $R$ in terms of $gr(\Gamma_E(R))$.

**Theorem 4.1:** Let $R$ be a finite commutative ring such that $gr(\Gamma_E(R)) = \infty$.
(i)   Then $Dim_d(\Gamma_E(R)) = 1$, if $R$ is a reduced ring,
(ii)  Then $Dim_d(\Gamma_E(R)) = 0$, if $R \cong \mathbb{Z}_4, \mathbb{Z}_9, \mathbb{Z}_2[x]/(x^2)$.

**Proof.** (a) Let us consider that $R$ is a reduced ring and not isomorphic to $\mathbb{Z}_2 \times \mathbb{Z}_2$, then it is certain that $R$ cannot be isomorphic to $\mathbb{Z}_2 \times A$, for some finite field $A$. Consequently, $R$ possesses two equivalence classes of ZD [(1,0)] & [(0,1)] which are adjacent to each other (by Theorem 2 [29]). As a result, $Dim_d(\Gamma_E(R)) =$

1. However, if $R$ is isomorphic to $\mathbb{Z}_2 \times \mathbb{Z}_2$, suggesting that $R$ is a Boolean ring, it follows that $\Gamma(R)$ is isomorphic to $\Gamma_E(R)$. Consequently, the result is completed by using (Lemma 3, [29]).

(b) For the rings listed as follows: $R \cong \mathbb{Z}_4, \mathbb{Z}_9, \mathbb{Z}_2[x]/(x^2)$, their CZDG, $\Gamma_E(R)$ represents a single vertex graph and hence $Dim_d(\Gamma_E(R)) = 0$. □

**Corollary 4.1:** Let $R$ be a finite commutative ring having unity 1 and $\mathbb{F}$ be a finite field. Also, let $R$ be a local ring isomorphic to any of the following listed rings, $\mathbb{Z}_8, \mathbb{Z}_{27}, \mathbb{Z}_2[x]/(x^3), \mathbb{Z}_4[x]/(2x, x^2 - 2), \mathbb{Z}_2[x,y]/(x^3, xy, y^2), \mathbb{Z}_8[x]/(2x, x^2), \mathbb{Z}_4[x]/(x^3, 2x^2, 2x), \mathbb{Z}_9[x]/(3x, x^2 - 6), \mathbb{Z}_9[x]/(3x, x^2 - 3), \mathbb{Z}_3[x]/(x^3)$. If $gr(\Gamma_E(R)) = \infty$, then for reduced rings $R \times \mathbb{F}$, $\Gamma_E(R \times \mathbb{F}) \cong \Gamma_E(R)$ with $Dim_d(\Gamma_E(R)) = 1$.

**Proof.** Let $gr(\Gamma_E(R)) = \infty$, then for all reduced rings $R \times \mathbb{F}$, we have $\Gamma_E(R \times \mathbb{F}) \cong k_{1,1}$ (Lemma 3, [29]). Furthermore, above local rings list has the same CZDG isomorphic to $k_{1,1}$. Hence $Mdim(\Gamma_E(R)) = 1$, for above local rings. □

We now focus on examining the connection between Ddim and diameter of CZDG. Considering, $diam(\Gamma_E(R)) \le 3$, when there is a cycle in CZDG, the following outcomes emerge.

**Theorem 4.2:** Let $R$ be a commutative ring and $\mathbb{F}_1$ and $\mathbb{F}_2$ are fields. Then,
(a) $Dim_d(\Gamma_E(R)) = diam(\Gamma_E(R)) = 1$ if $R \cong \mathbb{F}_1 \times \mathbb{F}_2$.
(b) $Dim_d(\Gamma_E(R)) = 0 \iff diam(\Gamma_E(R)) = 0$.
(c) $Dim_d(\Gamma_E(R)) = 0$ if $Z(R)^2 = 0$ and $|Z(R)| \ge 2$.
(d) $Dim_d(\Gamma_E(R)) = 0 \iff diam(\Gamma(R)) = 0$ or $1, R \ncong \mathbb{Z}_2 \times \mathbb{Z}_2$.

**Proof.**
(a) Let $R \cong \mathbb{F}_1 \times \mathbb{F}_2$, then $|\Gamma_E(R)| = 2$ (by Lemma 3, [29]). Since only equivalence classes of ZD are $[(0,1)]$ and $[(1,0)]$. So, CZDG is isomorphic to $k_{1,1}$. Hence, $Dim_d(\Gamma_E(R)) = 1 = diam(\Gamma_E(R))$.

(b) $Dim_d(\Gamma_E(R)) = 0$ iff $\Gamma_E(R)$ is single vertex graph and it is possible iff $diam(\Gamma_E(R)) = 0$.

(c) Let $|Z(R)| \ge 2$ and $Z(R)^2 = 0$. Hence $ann(a) = ann(b)$, for each $a, b \in Z(R)^*$, which implies that $diam(\Gamma_E(R)) = 0$. Therefore $Dim_d(\Gamma_E(R)) = 0$.

(d) Assume that $diam(\Gamma(R)) = 0$ or $1$, subsequently $\Gamma(R) \cong k_n$, thus $Dim_d(\Gamma_E(R)) = n - 1$ except if $R \ncong \mathbb{Z}_2 \times \mathbb{Z}_2$. Conversely, assume that $Dim_d(\Gamma_E(R)) = n - 1$. Then $\Gamma(R) \cong k_n$, hence $diam(\Gamma(R)) = 0$ or $1$. □

## 5. Conclusion

This research has effectively characterized rings by examining the dominant metric dimensions in their associated compressed zero-divisor graphs. The study establishes bounds for $Ddim(\Gamma_E(R))$ and thoroughly analyses Ddim for graphs $G$ that are realizable as $\Gamma_E(R)$. Furthermore, it explores the connections between dominant metric dimension, diameter, and girth of $\Gamma_E(R)$, enhancing our understanding of the intricate relationship between algebraic structures and graph representations. Based on the results, there is a promising avenue for future research to utilize Ddim in characterizing rings based on various types of zero-divisor graphs. The insights gained from this study serve as a solid foundation for further investigations, offering a platform for exploring different algebraic structures and their corresponding graph representations. The ongoing exploration of dominant metric dimensions holds the potential to uncover new perspectives and applications in mathematical research.

**Declaration:**

- Availability of data and materials: The data is provided on request to the authors.
- Conflicts of interest: The authors declare that they have no conflicts of interest and all the agree to publish this paper under academic ethics.
- Competing Interests and Funding: The authors did not receive support from any organization for the submitted work.
- Author's contribution: All the authors equally contributed towards this work.
- Acknowledgements: The author Nasir Ali, greatly appreciate the efforts of co-authors, the at most encouragement of his parents, Wife, and siblings.